# Explicit Construction of the Voronoi and Delaunay Cells of $W(A_n)$ and $W(D_n)$ Lattices and Their Facets

Mehmet Koca[a]
Nazife O. Koca[b], Abeer Al-Siyabi[c]

Department of Physics, College of Science, Sultan Qaboos University
P.O. Box 36, Al-Khoud, 123 Muscat, Sultanate of Oman
and
Ramazan Koc[d]
Department of Physics, Gaziantep University, 27310, Gaziantep, Turkey

**ABSTRACT**

Voronoi and Delaunay (Delone) cells of the root and weight lattices of the Coxeter-Weyl groups $W(a_n)$ and $W(d_n)$ are constructed. The face centered cubic (f.c.c.) and body centered cubic (b.c.c.) lattices are obtained in this context. Basic definitions are introduced such as parallelotope, fundamental simplex, contact polytope, root polytope, Voronoi cell, Delone cell, $n$-simplex, $n$-octahedron (cross polytope), $n$-cube and $n$-hemicube and their volumes are calculated. Voronoi cell of the root lattice is constructed as the dual of the root polytope which turns out to be the union of Delone cells. It is shown that the Delone cells centered at the origin of the root lattice $A_n$ are the polytopes of the fundamental weights $\omega_1, \omega_2, \ldots, \omega_n$ and the Delone cells of the root lattice $D_n$ are the polytopes obtained from the weights $\omega_1, \omega_{n-1}$ and $\omega_n$. A simple mechanism explains the tessellation of the root lattice by Delone cells. We prove that the $(n-1)$-facet of the Voronoi cell of the root lattice $A_n$ is $(n-1)$-dimensional rhombohedron and similarly the $(n-1)$-facet of the Voronoi cell of the root lattice $D_n$ is a dipyramid with a base of $(n-2)$-cube. Volume of the Voronoi cell is calculated via its $(n-1)$-facet which in turn can be obtained from the fundamental simplex. Tessellations of the root lattice with the Voronoi and Delone cells are explained by giving examples from lower dimensions. Similar considerations are also worked out for the weight lattices $A_n^*$ and $D_n^*$. It is pointed out that the projection of the higher dimensional root and weight lattices on the Coxeter plane leads to the $h$-fold aperiodic tiling where $h$ is the Coxeter number of the Coxeter-Weyl group. Tiles of the Coxeter plane can be obtained by projection of the two-dimensional faces of the Voronoi or Delone cells. Examples are given as the Penrose like 5-fold symmetric tessellation by the $A_4$ root lattice and the 8-fold symmetric tessellation by the $D_5$ root lattice.

[a]electronic-mail: mehmetkocaphysics@gmail.com ; retired professor
[b]electronic-mail: nazife@squ.edu.om
[c]electronic-mail: s21168@student.squ.edu.om
[d]electronic-mail: koc@gantep.edu.tr

## 1. Introduction

Higher dimensional crystallography described by the Coxeter-Dynkin diagrams [Engel, 1986; Engel, Michel & Senechal, 1994; Deza & Laurent, 1997] is not only a mathematical interest but also it is intimately related with quasi crystallography [de Brujin, 1981; Duneau & Katz, 1985; Baake, Joseph, Kramer & Schlottmann, 1990; Senechal, 1995; Chen, Moody & Patera, 1998] via projection into two and three dimensional Euclidean spaces [Koca, Koca & Koc, 2015; Boyle & Steinhardt, 2016].
Lie algebras derived from the root systems of the crystallographic groups are applied to particle physics as the symmetry of the standard model [Glashow, 1961; Weinberg, 1967; Salam, 1968; Fritzsch, Gell-Mann & Leutwyler, 1973] and its embeddings into higher-rank Lie groups [Georgi & Glashow, 1974; Fritzsch & Minkowski, 1975; Gursey, Ramond & Sikivie, 1976]. As a typical example, let us consider $SU(5)$ grand unified theory derived from $a_4$ root system that describes the unification of strong and electroweak interactions where the related polytopes describing the particle content project into 2-dimensions [Koca, Ozdes Koca & Koc, 2014]. As an application in quasicrystallography, in a recent paper [Koca, Ozdes Koca & Al-Siyabi, 2108], we have shown that the facet of the $A_4$ Voronoi cell is a rhombohedron and its 2- faces project into the Coxeter plane as *thick* and *thin* rhombuses of the Penrose tiling. The 2-faces of the Voronoi cell of the root lattice $D_5$ project on the Coxeter plane as two different triangles which lead to various 8-fold symmetric aperiodic tilings as we will demonstrate in Sec.4. This joint venture of higher dimensional crystallography and its symmetries invite further study of the Voronoi and Delone cells of the *A-D* lattices.

Voronoi cells [Voronoi, 1908 &1909], and Delaunay [Delaunay,1929 & 1938] polytopes have been studied in details in Chapter 21 of an excellent book by Conway and Sloane [Conway & Sloane, 1988] and especially in the paper [Conway and Sloane, 1991]. Further results on the Delone polytopes of the root lattices can be found in the reference [Deza & Grishukhin, 2004]. Numbers of facets of the Voronoi and Delone cells of the root lattice have been also determined by a technique of decorated Coxeter-Dynkin diagrams [Moody & Patera, 1992]. An expanded treatment of the lattices derived from Coxeter-Weyl groups can be found in [Engel, Michel & Senechal, 1994]. However detailed structures of the $(n-1)$-facets of the Voronoi cells have not been studied in details whose 2-faces are essential for the aperiodic tiling of the Coxeter plane which could be taken as models for the quasicrystallography. One may also find detailed discussions of the polytopes in two major references [Coxeter, 1973] and [Grunbaum, 1967]. Lie algebraic technique of the lattices derived from the Coxeter-Dynkin diagrams have been studied in [Bourbaki, 1968] and [Humphreys, 1992].

In the present paper, we construct the Voronoi cell $V(0)$ centered around the origin, as the dual polytope of the root polytope determined by the root system of the associated Lie algebra. Being dual of each other, the number of $d$-dimensional facet $N_d^r$ of the root polytope ($r$ stands for the root polytope) equals the number of $(n-d-1)$-dimensional facet $N_{n-d-1}^v$ of the Voronoi cell ($v$ stands for the Voronoi polytope) which follows directly from the Coxeter-Dynkin diagram. It turns out that the Voronoi cell $V(0)$ is disjoint union of Delone cells. For the root lattice $A_n$ Voronoi cell is the union of orbits of the fundamental weights $\omega_1, \omega_2, ..., \omega_n$ and the Voronoi cell of the root lattice $D_n$ consists of the orbits of the weights $\omega_1, \omega_{n-1}$ and $\omega_n$.

The paper is organized as follows. In Sec.2, we introduce the main tools of the lattices such as root system, weight vectors, parallelotope, fundamental simplex, contact polytope, Voronoi cell and Delone cells derived from the Coxeter-Dynkin diagrams and important aspects of the associated Coxeter-Weyl groups are discussed. Sec.3 deals with the root lattice $A_n$ by introducing the root polytope and its dual. We prove that the $(n-1)$-dimensional facets of the Voronoi cell are generalized rhombohedra and give examples of dimensions $n = 1, 2, 4$. We calculate the volume of the Voronoi cell via its facets. We study in Sec.4 the same problem for the Voronoi cell of the root lattice $D_n$ and show that $(n-1)$-dimensional facet of the Voronoi cell is a dipyramid with a base of $(n-2)$-dimensional cube. The lattices $D_3$, $D_4$ and $D_5$ are studied as typical examples. In Sec.5 we study the weight lattices $A_n^*$ and $D_n^*$ and construct their Voronoi and Delone cells and contact polytopes. Appendix A involves derivations of the formula for the volumes of the $n$-simplex $\alpha_n$, $n$-octahedron $\beta_n$ (cross polytope) and $n$-dimensional hemicube $h\gamma_n$.

## 2. Lattices and Polytopes derived from Coxeter-Dynkin diagrams

This paper will not discuss the lattices of the exceptional Lie algebras that deserve a separate exposition as they are also related to octonionic and quaternionic representations [Coxeter, 1946; Koca & Ozdes, 1989; Koca, Koc & Al-Barwani, 2006; Koca, 2007]. Neither we will consider the lattices generated by the short roots of the root system of the B-C series for the $B_n$ lattice is the simple cubic lattice and the $C_n$ can be represented by the $D_n$ lattice [Conway & Sloane, 1988]. The Coxeter-Dynkin diagrams and their extended diagrams of the A-D series are shown in Fig.1.

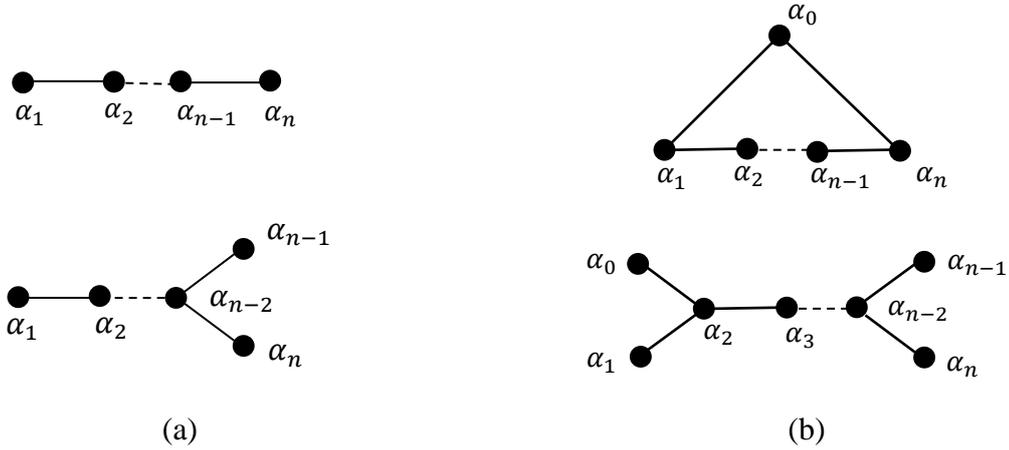

**Figure 1**
(a) Coxeter-Dynkin diagrams of $a_n$ and $d_n$, (b) Extended Coxeter Dynkin diagrams of $a_n$ and $d_n$.

The nodes represent the simple roots $\alpha_i, (i = 1, 2, ..., n)$ of the associated Lie algebra of rank $n$ where the norm of the roots are given by $(\alpha_i, \alpha_i) = 2$ and we define the Cartan matrix (Gram matrix in lattice terminology) by the relation

$$C_{ij} = \frac{2(\alpha_i, \alpha_j)}{(\alpha_j, \alpha_j)}. \qquad (1)$$

The fundamental weight vectors $\omega_i$ are defined by the relation $(\omega_i, \alpha_j) = \delta_{ij}$ where $\delta_{ij}$ is the Kronecker-delta and they are related to each other by the relations

$$\alpha_i = \sum_j C_{ij} \omega_j, \quad \omega_i = \sum_j (C^{-1})_{ij} \alpha_j, \qquad (2)$$

where the scalar product of the fundamental weights define the matrix elements of the inverse Cartan matrix $(\omega_i, \omega_j) = (C^{-1})_{ij}$.

The root lattice $\Lambda$ is defined as the set of vectors $p = \sum_{i=1}^{n} b_i \alpha_i$, $b_i \in \mathbb{Z}$. Among many other tessellations the lattice is tiled with the *parallelotope* generated by the simple roots $\alpha_i$ which, in an orthogonal base, constitute the rows of the generator matrix $M$ where the Cartan matrix can be written as $C = MM^T$. This implies that the volume of the parallelotope is $\sqrt{detC}$ which is also called the volume of the lattice. As we will see later that $\sqrt{detC}$ is the volume of the Voronoi cell. A similar relation applies to the weight lattice $G^*$ whose vectors are the linear combinations of the weight vectors $q = \sum_{i=1}^{n} c_i \omega_i$, $c_i \in \mathbb{Z}$ and volume of the parallelotope is $1/\sqrt{detC}$. Let $r_i$, $(i = 1, 2, ..., n)$ denote the reflection generator with respect to the hyperplane orthogonal to the simple root $\alpha_i$ which operates on an arbitrary vector $\lambda$ as

$$r_i \lambda = \lambda - \frac{2(\lambda, \alpha_i)}{(\alpha_i, \alpha_i)} \alpha_i. \tag{3}$$

It transforms a fundamental weight vector as $r_i \omega_j = \omega_j - \alpha_i \delta_{ij}$. The reflection generators generate the Coxeter group $W(g) = < r_1, r_2, ..., r_n | (r_i r_j)^{m_{ij}} >$ which is also called the Coxeter-Weyl group for the crystallographic Coxeter groups. Now we will introduce the definitions of the basic ingredients of the lattice terminology.

*Root polytope* (*root system*): The point group $W(g)$ acting on one of the simple roots of a simly-laced Lie algebra generate the root system $W(g)\alpha_i$. *The root polytope* is a convex polytope whose vertices are vectors of the root system. The root polytope of the root lattice $A_n$ is generated also from the highest weight of the adjoint representation of the Lie algebra $a_n$ by $W(a_n)(\omega_1 + \omega_n)$. The root polytope of the group $W(d_n)$ is the orbit $W(d_n)(\omega_2)$. The highest weight vector [Slansky,1981] for an irreducible representation of the Lie algebra is defined as the weight vector $q = \sum_{i=1}^{n} c_i \omega_i =: (c_1, c_2, ..., c_n)$, ($c_i$ integer with $c_i \geq 0$).

We denote by the convex polytope possessing the symmetry of the Coxeter-Weyl group as the orbit of the highest weight vector

$$W(g)(c_1, c_2, ..., c_n) =: (c_1, c_2, ..., c_n)_g. \tag{4}$$

In the Coxeter-Dynkin diagram, we will use the component $c_i$ of $\omega_i$ on top of the corresponding node as shown in the diagrammatic notations of Fig. 2.

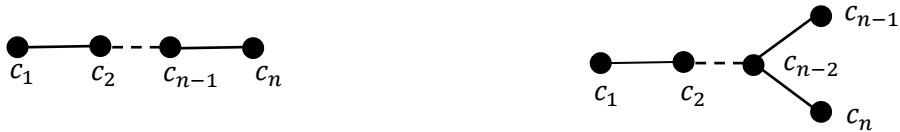

**Figure 2**
Diagrammatic notation for an arbitrary convex polytope of $a_n$ and $d_n$ respectively.

With this notation the root polytope of $W(a_n)$ will be shown either by $(10...01)_{a_n}$ or simply by $(\omega_1 + \omega_n)_{a_n}$ or we use the diagrammatic notation of Fig. 2. Similarly the root polytope of $W(d_n)$ is denoted either by $(010...00)_{d_n}$ or by $(\omega_2)_{d_n}$.

Another polytope associated with a lattice is the *contact polytope*. It is the root polytope for the root lattice but it has more general definition: consider a sphere packing of a lattice and the centers of the spheres touching to the central sphere. *Contact polytope* is the convex hull of the centers of these spheres. We will identify the contact polytope in the weight lattices $A_n^*$ and $D_n^*$.

The *Affine Weyl group* $W(G)$ is the infinite discrete group generated by $W(g)$ and an additional generator, usually denoted by $r_0$, which describes the reflection with respect to the hyperplane bisecting the highest weight vector: $W(G) = <r_0, r_1, r_2, ..., r_n>$. The generator $r_0$ acts as a translation by taking zero vector (0) to the highest weight of the adjoint representation. The point group $W(g)$ is an invariant subgroup of the affine Weyl group $W(G)$ where the root lattice is denoted by the letter $G$.

The *fundamental simplex* of the lattice $A_n$ is a convex polytope with $n+1$ vertices given by $\omega_1, \omega_2, ..., \omega_n$ and the origin (0). *The roof of the fundamental simplex* is the polytope determined by the vertices $\omega_1, \omega_2, ..., \omega_n$ excluding the origin. As we will see later the Voronoi cell of the lattice $A_n$ is the union of the orbits

$$(\omega_1)_{a_n} \cup (\omega_2)_{a_n} ... \cup (\omega_n)_{a_n}. \qquad (5)$$

We will discuss the structure of the fundamental simplex of the lattice $D_n$ in Sec. 4.

The *Voronoi polytope* $V(p)$ centered at a lattice point $p \in \Lambda$ is the set of points

$$V(p) = \{x \in \mathbb{R}^n : \| x - p \| \leq \| x - q \|, \text{for all } q \in \Lambda\}, \qquad (6)$$

where $\| x \| = \sqrt{(x,x)}$ is the length of $x$. The Voronoi polytopes are all congruent so we will study only $V(0)$, the polytope centered at the origin which can also be defined as

$$V(0) = \left\{ x \in \mathbb{R}^n, \forall p, (x,p) \leq \tfrac{1}{2}(p,p) \right\}. \qquad (7)$$

This definition shows that the facet of the Voronoi cell is the hyperplane bisecting the orthogonal vector $0p$.

Let $v \in \mathbb{R}^n$ be a vertex of an arbitrary Voronoi polytope $V(p)$. The convex hull of the lattice points closest to $v$ is called the *Delone polytope* containing $v$. Vertices of the Voronoi cell $V(0)$ of the root lattices of the $A_n$-$D_n$ series consist of the vertices of the Delone cells centered at the origin. There is an interesting relation between the Voronoi cell $V(0)$ and the Delone cells, namely the volume of $V(0)$ equals the sum of the volumes of the Delone cells centered at the origin [ Michel, 1994].

Some regular polytopes will be of special interest like *n-simplex $\alpha_n$*, *cross polytope $\beta_n$*, (*n-octahedron*), *n-hemicube $h\gamma_n$* and *n-cube $\gamma_n$*. They are represented as the polytopes $(\omega_1)_{a_n}$ (or $(\omega_n)_{a_n}$, $(\omega_1)_{d_n}$, $(\omega_n)_{d_n}$ (or $(\omega_{n-1})_{d_n}$) and the union $(\omega_n)_{d_n} \cup (\omega_{n-1})_{d_n}$ respectively. The volumes of these polytopes are calculated in Appendix A by using recurrence relations.

## 3. The Root Lattice $A_n$ and related polytopes

A useful representation of the roots and the weights of $A_n$ can be given in terms of an orthonormal set of vectors, $l_i$ $(i = 1, 2, ..., n + 1)$, $(l_i, l_j) = \delta_{ij}$. We define the simple roots as linear combinations of orthonormal vectors $\alpha_i = l_i - l_{i+1}$, $(i = 1, 2, ..., n)$. The group generators $r_i$ permute the set of orthonormal vectors as $r_i: l_i \leftrightarrow l_{i+1}$. Clearly, the Coxeter-Weyl group $W(a_n) \simeq S_{n+1}$ represents the permutation group of the $n + 1$ objects with an order of $(n + 1)!$. When we define the vectors in terms of their components in the $n + 1$ dimensional Euclidean space the simple roots and the fundamental weights read

$$\alpha_1 = (1, -1, 0, ..., 0); \quad \alpha_2 = (0, 1, -1, 0, ..., 0); \quad ...; \quad \alpha_n = (0, 0, ..., 1, -1);$$
$$\omega_1 = \frac{1}{n+1}(n, -1, -1, ..., -1); \quad \omega_2 = \frac{1}{n+1}(n-1, n-1, -2, -2, ..., -2); \quad ...; \quad (8)$$
$$\omega_n = \frac{1}{n+1}(1, 1, ..., 1, -n).$$

In an elegant notation the fundamental weights read $\omega_i = \frac{1}{n+1}\left((j)^i, (-i)^j\right)$, where $i + j = n + 1$ [Conway& Sloane, 1988]. These are called holes of the root lattice $A_n$ denoted by a much simpler notation $[i] \equiv \omega_i$. To illuminate the topic, first we will give some examples.

*i)* The lattice $A_1 \approx A_1^*$

The simplest lattice of course is the lattice $A_1 \approx A_1^*$ with a point group of order 2 generated by $r$, the reflection with respect to the origin. Here the root polytope is the line segment between the points $\alpha$ and $-\alpha$. Since the weight is $\omega = \frac{\alpha}{2}$ the Voronoi cell is the line segment between $\omega$ and $-\omega$. The generator $r_o$ reflects the points in the straight line with respect to the point $\omega = \frac{\alpha}{2}$ corresponding to the hyperplane bisecting the highest weight $\alpha$. The Delone polytopes are the line segments between 0 and $\pm \alpha$ and those obtained by translation by the integer multiples of the simple root $\alpha$. Tessellation by the Voronoi cell is the translation of the Voronoi interval by integers. If $\alpha$ is scaled by $\sqrt{2}$ the weight lattice corresponds to the eigenvalues of the $S_z$ operator of the $SU(2)$ Lie algebra with values $..., -\frac{3}{2}, -1, -\frac{1}{2}, 0, \frac{1}{2}, 1, \frac{3}{2}, ...$ and the root lattice corresponds to the eigenvalues of the integer values of $S_z$ $... - 3, -2, -1, 0, 1, 2, 3 ...$. Although the root lattice $A_1$ is a sublattice of the weight lattice $A_1^*$, they are equivalent in the sense that the root lattice can be obtained from the weight lattice by multiplying the weight lattice by the factor 2. The fundamental simplex is the line segment between 0 and $\omega$ where the single point $\omega$ represents the roof of the simplex.

*ii)* $A_2 \approx A_2^*$

Let us also look at the next relatively simpler case $A_2$. The hexagonal lattice $A_2$ is quite well known as it represents the graphene [Novoselov et al, 2004]. The root polytope is a regular hexagon with 6 vertices and 6 edges obtained as the orbit of the weight vector $\omega_1 + \omega_2 = \alpha_1 + \alpha_2$ as shown in Fig. 3. A typical face is the line segment determined by the points $\omega_1 + \omega_2$ and $r_1(\omega_1 + \omega_2) = \omega_1 + \omega_2 - \alpha_1$. Dual vector to the line segment $\alpha_1$ is the weight vector $\omega_2$, the orbit of which is the triangle $(\omega_2)_{a_2} = \{\omega_2, \omega_1 - \omega_2, -\omega_1\}$. Similarly the vector orthogonal to the line segment $\alpha_2$ is $\omega_1$ from which one generates the triangle $(\omega_1)_{a_2} = -(\omega_2)_{a_2}$. The dual polytope of the root polytope is then the regular hexagon which is disjoint union of two orbits representing the Voronoi cell $V(0) = (\omega_1)_{a_2} \cup (\omega_2)_{a_2}$ as depicted in Fig. 3. The triangle $(\omega_1)_{a_2}$ and the inverted triangle $(\omega_2)_{a_2}$ each represents a Delone cell centered at the origin. Area of each triangle is $\frac{\sqrt{3}}{2}$ so the sum of two areas $\sqrt{3} = \sqrt{\det C_2}$ equals the volume (area) of the Voronoi cell of $A_2$. The fundamental simplex is the triangle (Fig. 3 (b)) (dashed lines) with the vertices 0, $\omega_1$ and $\omega_2$ where the roof is the line segment between $\omega_1$ and $\omega_2$. The vertices of the Delone cells centered at the vertices of the Voronoi cell $V(0)$ are obtained by adding the vertices of two sets of the triangles to each other: $\{\omega_2, \omega_1 - \omega_2, -\omega_1\} + (-\omega_2, -\omega_1 + \omega_2, \omega_1)$. Then the lattice $A_2$ is generated by translation. As it is clear from Fig. 3 the weight lattice $A_2^*$ generated by the Voronoi cell $V(0)$ will be a scaled copy of the root lattice. Hence, two lattices are congruent. The Voronoi cell of $A_2^*$ is the hexagon represented by the orbit $V(0)^* = \frac{1}{3}(\omega_1 + \omega_2)_{a_2}$.

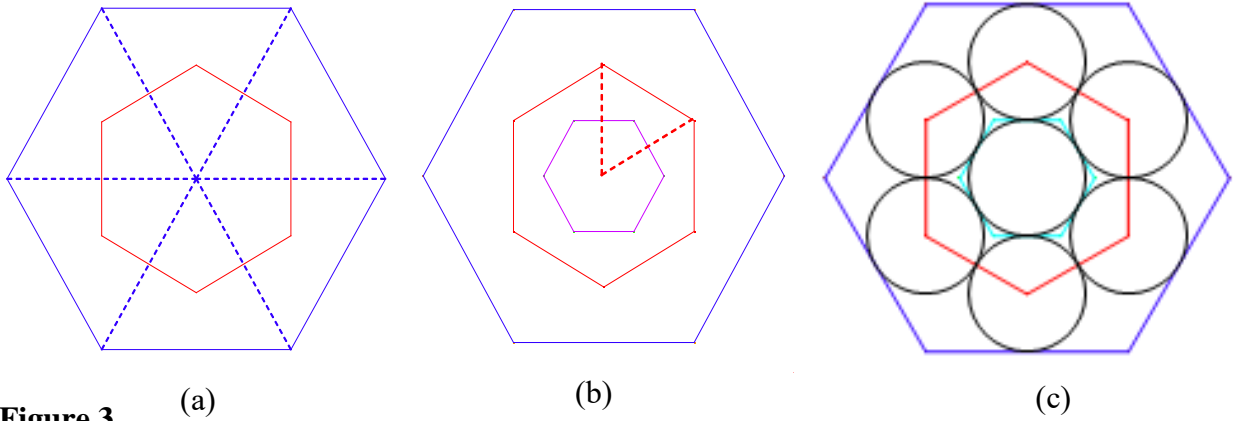

**Figure 3**
(a) Root polytope, Voronoi cell and Delone cells of $A_2$ lattice, (b) Root polytope, Voronoi cell of $A_2$, Voronoi cell of $A_2^*$ and fundamental simplex, (c) Sphere packing and related polytopes.

*iii*) $A_4$

We will discuss the lattice $A_3 \approx D_3$ in Sec.4 as it is the crossroads between two lattices $A_n$ and $D_n$. The Coxeter-Dynkin diagram of the group $W(a_4)$ can be inferred from Fig.1 and Fig. 2. The order of the group is $|W(a_4)| = 120$ which can be extended to the automorphism group of order 240 by the generator $\gamma$: $\omega_1 \leftrightarrow \omega_4$, $\omega_2 \leftrightarrow \omega_3$. Number of a particular facet of a polytope is determined as the index of the subgroup of the Coxeter-Weyl group leaving the facet invariant. The highest weight $\omega_1 + \omega_4$ is left invariant under the subgroup $W(a_2) = <r_2, r_3>$. The root polytope is then 4-dimensional convex polytope with numbers of facets denoted respectively by $N_0$ vertices, $N_1$ edges, $N_2$ polygons and $N_3$ polyhedra given by,

$$N_0^r = \frac{|<r_1,r_2,r_3,r_4>|}{|<r_2,r_3>|} = 20,$$

$$N_1^r = \frac{|<r_1,r_2,r_3,r_4>|}{|<r_1><r_3>|} + \frac{|<r_1,r_2,r_3,r_4>|}{|<r_4><r_2>|} = 60,$$

$$N_2^r = \frac{|<r_1,r_2,r_3,r_4>|}{|<r_1,r_2>|} + \frac{|<r_1,r_2,r_3,r_4>|}{|<r_3,r_4>|} + \frac{|<r_1,r_2,r_3,r_4>|}{|<r_1><r_4>|} \qquad (9)$$
$= 20\text{ (triangles)} + 20\text{(triangles)} + 30\text{ (squares)} = 70,$

$$N_3^r = \frac{|<r_1,r_2,r_3,r_4>|}{|<r_1,r_2,r_3>|} + \frac{|<r_1,r_2,r_3,r_4>|}{|<r_2,r_3,r_4>|} + \frac{|<r_1,r_2,r_3,r_4>|}{|<r_1,r_2><r_4>|} + \frac{|<r_1,r_2,r_3,r_4>|}{|<r_3,r_4><r_1>|}$$
$= 10\text{(tetrahedra)} + 20\text{(triangular prisms)} = 30.$

This detailed formulation is given for the convenience of the reader. Now the facets of the Voronoi cell $V(0)$ equals $N_0^v = 30$ vertices, $N_1^v = 70$ edges, $N_2^v = 60$ rhombuses and $N_3^v = 20$ rhombohedra. These numbers satisfy the Euler characteristic equation $N_0 - N_1 + N_2 - N_3 = 0$.

To determine the vertices of the dual polytope we note that the group $<r_1, r_2, r_3>$ generating a tetrahedron acting on $\omega_1 + \omega_4$ leaves the vector $\omega_4$ invariant. Similarly the groups $<r_2, r_3, r_4>, <r_1, r_2><r_4>$ and $<r_3, r_4><r_1>$ leave the vectors $\omega_1, \omega_3, \omega_2$ invariant respectively. It is straight forward to show that the hyperplane determined by the fundamental weights is orthogonal to the highest weight vector $\omega_1 + \omega_4$. The union of the orbits of these weights constitute 30 vertices of the Voronoi cell $V(0) = (\omega_1)_{a_4} \cup (\omega_2)_{a_4} \cup (\omega_3)_{a_4} \cup (\omega_4)_{a_4}$.

Volume of the Voronoi cell can be obtained as

$$\text{Vol } V(0) = \sum_{i=1}^{4} Vol\ (\omega_i)_{a_4} = \frac{\sqrt{5}}{24} + \frac{11\sqrt{5}}{24} + \frac{11\sqrt{5}}{24} + \frac{\sqrt{5}}{24} = \sqrt{5}. \tag{10}$$

Each orbit constituting the Voronoi cell $V(0)$ represents a Delone cell centered at the origin: $(\omega_1)_{a_4}$ (4-simplex), $(\omega_2)_{a_4}$ (ambo 4-simplex), $(\omega_3)_{a_4} = -(\omega_2)_{a_4}$ ($2^{nd}$ ambo 4-simplex), $(\omega_4)_{a_4} = -(\omega_1)_{a_4}$ ($3^{rd}$ ambo 4-simplex). See for the definitions [Conway & Sloane, 1991]. The 4-simplex sometimes called the 5-cell to remind that it consists of 5 tetrahedra as facets. The vertices of the ambo 4-simplex $(\omega_2)_{a_4}$ are the pairwise sum of the vertices of the 4-simplex $(\omega_1)_{a_4} = \{\omega_1, \omega_2 - \omega_1, \omega_3 - \omega_2, \omega_4 - \omega_3, -\omega_4\}$ and it consists of tetrahedra and octahedra as facets. The others follow from the above relations. Tessellation of the root lattice with Delone cells can be simply explained as follows. For example, when the vertices of the polytope $(\omega_4)_{a_4}$ are added to the vertices of the 4-simplex $(\omega_1)_{a_4}$ one obtains 5 Delone cells whose vertices are at the root lattice and centered at the vertices of the polytope $(\omega_4)_{a_4}$ of the Voronoi cell. Similarly, the vertices of the Delone cell centered at the vertices of $(\omega_3)_{a_4}$ can be obtained by adding its vertices to the vertices of Delone cells of the polytope $(\omega_2)_{a_4}$. With the negatives of the above Delone cells we obtain a tiling of the root lattice with Delone cells centered at the vertices of the Voronoi cell $V(0)$. The complete tessellation is carried out by the affine Coxeter-Weyl group $W(A_n)$.

Now we discuss one of the facet of the Voronoi cell $V(0)$. The vertices of the 3-facet of $V(0)$ centered at the vertex $\frac{1}{2}(\omega_1 + \omega_4)$ can be obtained by the subgroup $<r_2, r_3>$ operating on the set of weights $\{\omega_1, \omega_2, \omega_3, \omega_4\}$. Since the group $<r_2, r_3>$ leaves the weights $\omega_1$ and $\omega_4$ invariant we obtain a polyhedron with 8 vertices as

$<r_2, r_3> \{\omega_1, \omega_2, \omega_3, \omega_4\} = \{\omega_1, \omega_2, \omega_3, \omega_4, (r_2 r_3)\omega_2, (r_2 r_3)\omega_3, (r_3 r_2)\omega_2, (r_3 r_2)\omega_3\}.$ (11)

The polyhedron is a rhombohedron generated by 3 edges $\{\omega_1 - \omega_2, \omega_2 - \omega_3, \omega_3 - \omega_4\}$ as shown in Fig. 4 and its volume is $\sqrt{\frac{2}{5}}$.

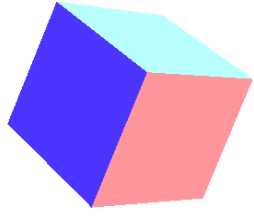

**Figure 4**
The rhombohedron as a facet of the Voronoi polytope $V(0)$ of $A_4$.

One can easily check that generating three edges are equal

$$\| \omega_1 - \omega_2 \| = \| \omega_2 - \omega_3 \| = \| \omega_3 - \omega_4 \| = \frac{2}{\sqrt{5}} \tag{12}$$

and the angle between any pair is $\phi = cos^{-1}(-\frac{1}{4}) \approx (104.5)^0$. All 2-faces of the Voronoi cell are rhombuses and the 3-facets are rhombohedra. It was pointed out that 2-faces project onto the Coxeter plane as thin and tick rhombuses of the Penrose tiling [Koca, Ozdes Koca & Al-Siyabi, 2018].

Volume of the Voronoi cell $V(0)$ can also be computed as the sum of the volumes of four-dimensional pyramids based on 20 rhombohedra as

$$20 \left(\frac{1}{4}\right) \left(\frac{\|\omega_1 + \omega_4\|}{2}\right) \sqrt{\frac{2}{5}} = \sqrt{5} . \tag{13}$$

The roof of the fundamental simplex is a semiregular tetrahedron consisting of 2-faces of isosceles triangles of edges $(\frac{2}{\sqrt{5}}, \frac{2}{\sqrt{5}}, \sqrt{\frac{6}{5}})$ and 2 other faces of isosceles triangles of edges $(\frac{2}{\sqrt{5}}, \sqrt{\frac{6}{5}}, \sqrt{\frac{6}{5}})$ and its volume is $\frac{1}{6}\sqrt{\frac{2}{5}}$.

The fundamental simplex consists of 5-cells identical to the roof of the semiregular tetrahedron. We also note that the volume of the fundamental simplex is $\left(\frac{1}{4}\right) \left(\frac{\|\omega_1 + \omega_4\|}{2}\right) \frac{1}{6}\sqrt{\frac{2}{5}} = \frac{1}{4!\sqrt{5}}$.

It is evident that the volume of the Voronoi cell is $Vol\ V(0) = 5! \times Vol$ (fundamental simplex).

After these basic examples we can discuss the lattice $A_n$ in its general context.

The root system representing the vertices of the root polytope can be obtained from the diagram of Fig.1 as $l_i - l_j$, $(i,j = 1, 2, ..., n+1)$ and the number of vertices equals $N_0^r = \frac{(n+1)!}{(n-1)!} = n(n+1)$. Similarly number of edges is given by $N_1^r = \frac{(n+1)!}{(n-2)!}$. It can be shown that the number of $d$-facets is given by

$$N_d^r = \frac{(n+1)!}{(n-1-d)!(d+2)!}(2^{2+d} - 2). \tag{14}$$

One can reproduce the values in (9) by substituting $n = 4$ in (14). Since the facets of the dual polytope is given by $N_d^v = N_{(n-d-1)}^r$ the number of facets of the Voronoi cell is obtained as

$$N_d^v = \frac{(n+1)!}{(n+1-d)!(d)!}(2^{n+1-d} - 2). \tag{15}$$

Compare with the result of the reference [Moody & Patera, 1992].

We will now prove that all 2-faces of the Voronoi cell is a rhombus and any $d$-facet for $3 \leq d \leq n-1$ is a $d$-dimensional rhombohedron. Let us recall that the hyperplane orthogonal to the vector $\frac{(\omega_1 + \omega_n)}{2}$ is given by the set of fundamental weights $\{\omega_1, \omega_2, \ldots, \omega_n\}$ for $((\omega_i - \omega_j), (\omega_1 + \omega_n)) = 0$. The vertices of the $(n-1)$-facet centered at $\frac{(\omega_1 + \omega_n)}{2}$ is given by the orbit

$$< r_2, r_3, \ldots, r_{n-1} > \{\omega_1, \omega_2, \ldots, \omega_{n-1}, \omega_n\}. \tag{16}$$

We will prove that this is a rhombohedron in $(n-1)$-dimensional Euclidean space. We list the number of vertices in Table.1 generated from each fundamental weight.

**Table 1**
Number of vertices of $(n-1)$-facet generated from fundamental weights.

| Fundamental weights | Number of generated vertices |
|---|---|
| $\omega_1$ | $\binom{n-1}{0}$ |
| $\omega_2$ | $\binom{n-1}{1}$ |
| $\omega_3$ | $\binom{n-1}{2}$ |
| $\vdots$ | $\vdots$ |
| $\omega_{n-1}$ | $\binom{n-1}{n-2}$ |
| $\omega_n$ | $\binom{n-1}{n-1}$ |

Total number of vertices is given by $\sum_{i=0}^{n-1}\binom{n-1}{i} = 2^{n-1}$.

The $(n-1)$-dimensional rhombohedron can be generated by the $(n-1)$ vectors given by

$$k_1 = \omega_1 - \omega_2, \; k_2 = \omega_2 - \omega_3, \ldots, k_{n-1} = \omega_{n-1} - \omega_n \tag{17}$$

having the same length $\| k_i \| = \sqrt{\frac{n}{n+1}}, (i = 1, 2, \ldots, n-1)$ and the equal angle between any pair of generating vectors

$$(k_i, k_j) = -\frac{1}{n+1}, i \neq j, \quad \phi = cos^{-1}(-\frac{1}{n}). \tag{18}$$

This proves that the two-dimensional face is a rhombus, three-dimensional face is a rhombohedron and the higher dimensional polyhedron is the higher dimensional rhombohedron. The vectors $k_i$ have $n+1$ components in the basis of the vectors $l_i$. However it is possible to introduce a new set of orthonormal vectors

$$l_1' = \frac{l_1 - l_{n+1}}{\sqrt{2}}, \quad l_2' = \frac{l_1 + l_{n+1} - 2l_2}{\sqrt{6}}, \ldots,$$
$$l_{n-1}' = \frac{l_1 + l_{n+1} + l_2 + \cdots - (n-1)l_{n-1}}{\sqrt{n(n-1)}}, \quad l_n' = \frac{l_1 + l_{n+1} + l_2 + \cdots - nl_n}{\sqrt{(n+1)n}}, \quad l_{n+1}' = \frac{l_1 + l_{n+1} + l_2 + \cdots + l_n}{\sqrt{n+1}} \tag{19}$$

where the vectors $k_i$ are expressed in terms of the linearly independent $(n-1)$ components of $l_i', i = 2, 3, \ldots, n$. Then the generator matrix

$$M = \begin{bmatrix} k_{11} & k_{12} & \cdots & k_{1m} \\ \vdots & \ddots & & \vdots \\ k_{m1} & k_{m2} & \cdots & k_{mm} \end{bmatrix} \tag{20}$$

is given by the $(n-1) \times (n-1)$ matrix as

$$M_{n-1} = \begin{pmatrix} \sqrt{2/3} & -\frac{1}{\sqrt{12}} & -\frac{1}{\sqrt{20}} & -\frac{1}{\sqrt{30}} & -\frac{1}{\sqrt{42}} & \cdots & -\frac{1}{\sqrt{(n-1)^2+(n-1)}} & -\frac{1}{\sqrt{n^2+n}} \\ 0 & \sqrt{3/4} & -\frac{1}{\sqrt{20}} & -\frac{1}{\sqrt{30}} & -\frac{1}{\sqrt{42}} & \cdots & -\frac{1}{\sqrt{(n-1)^2+(n-1)}} & -\frac{1}{\sqrt{n^2+n}} \\ 0 & 0 & \sqrt{4/5} & -\frac{1}{\sqrt{30}} & -\frac{1}{\sqrt{42}} & \cdots & -\frac{1}{\sqrt{(n-1)^2+(n-1)}} & -\frac{1}{\sqrt{n^2+n}} \\ 0 & 0 & 0 & \sqrt{5/6} & -\frac{1}{\sqrt{42}} & \cdots & -\frac{1}{\sqrt{(n-1)^2+(n-1)}} & -\frac{1}{\sqrt{n^2+n}} \\ 0 & 0 & 0 & 0 & \sqrt{6/7} & \cdots & -\frac{1}{\sqrt{(n-1)^2+(n-1)}} & -\frac{1}{\sqrt{n^2+n}} \\ . & . & . & . & . & \cdots & . & . \\ . & . & . & . & . & \cdots & . & . \\ 0 & 0 & 0 & 0 & 0 & \cdots & \sqrt{\frac{n-1}{n}} & -\frac{1}{\sqrt{n^2+n}} \\ 0 & 0 & 0 & 0 & 0 & \cdots & 0 & \sqrt{\frac{n}{n+1}} \end{pmatrix} \tag{21}$$

Clearly, the det $M = \sqrt{\frac{2}{n+1}}$ is the volume of the $(n-1)$ dimensional rhombohedron.
Volume of the Voronoi cell $V(0)$ can then be determined as

$$Vol\, V(0) = n(n+1) \times \frac{1}{n} \times \left(\frac{\|\omega_1 + \omega_n\|}{2}\right)\sqrt{\frac{2}{n+1}} = \sqrt{n+1}. \tag{22}$$

Since the rhombohedron is obtained from (16) it can be partitioned into $(n-1)!$ identical polytopes congruent to the roof of the fundamental simplex. Therefore, the volume of the roof of the fundamental simplex equals $\frac{1}{(n-1)!}\sqrt{\frac{2}{n+1}}$ and with a similar formula to (22) the volume of the fundamental simplex is given by

$$Vol(\text{fundamental simplex}) = \frac{1}{n} \times \left(\frac{\|\omega_1 + \omega_n\|}{2}\right)\sqrt{\frac{2}{n+1}} \times \frac{1}{(n-1)!} = \frac{1}{n!\sqrt{n+1}}. \tag{23}$$

It is then obvious that $Vol\ V(0) = (n+1)!\ Vol($ fundamental simplex$)$.

In regard to the Delone cells we recall that the orbits $(\omega_1)_{a_n}, (\omega_2)_{a_n}, \ldots, (\omega_n)_{a_n}$ constitute the Delone cells centered at the origin. Tessellation of the root lattice $A_n$ with Delone cells can be carried out by adding the vectors of the pairs of orbits $(\omega_1)_{a_n} + (\omega_n)_{a_n}, (\omega_2)_{a_n} + (\omega_{n-1})_{a_n}, \ldots, (\omega_{n/2})_{a_n} + (\omega_{n+1-\frac{n}{2}})_{a_n}$ for even $n$. For odd $n$ the vectors of the orbit $(\omega_{\frac{n+1}{2}})_{a_n}$ will be added to themselves. By this technique not only we obtain the Delone cells with the vertices of the root lattice but also their centers are the vertices of the Voronoi cell $V(0)$. Applying the affine Coxeter-Weyl group $W(A_n)$ on the Delone cells we obtain the $A_n$ lattice tessellated by the Delone cells or applying the same group on the vertices of the Voronoi cell the root lattice $A_n$ will be tessellated with the Voronoi cell.

## 4. The Root Lattice $D_n$ and related polytopes

The simple roots and the fundamental weights are given as follows

$$\alpha_1 = l_1 - l_2,\ \alpha_2 = l_2 - l_3, \ldots, \alpha_{n-1} = l_{n-1} - l_n,\ \alpha_n = l_{n-1} + l_n,$$

$$\omega_1 = l_1,\ \omega_2 = l_1 + l_2,\ \omega_3 = l_1 + l_2 + l_3, \ldots,\ \omega_{n-1} = \tfrac{1}{2}(l_1 + l_2 + \cdots + l_{n-1} - l_n), \quad (24)$$

$$\omega_n = \tfrac{1}{2}(l_1 + l_2 + \cdots + l_{n-1} + l_n).$$

Generators of the Coxeter-Weyl group transform the orthonormal set of vectors as $r_i: l_i \leftrightarrow l_{i+1}, (i = 1, 2, \ldots, n-1)$ and $r_n: l_{n-1} \leftrightarrow -l_n$. Order of the point group is $|W(d_n)| = 2^{n-1}n!$. The root polytope is the orbit $W(d_n)(\omega_2) \equiv (\omega_2)_{d_n}$ which consists of $2n(n-1)$ vectors,

$$(\omega_2)_{d_n} = \{\pm l_i \pm l_j\}, (i \neq j = 1, 2, \ldots, n). \quad (25)$$

The root lattice $D_n$ consists of the vectors $\sum_{i=1}^n m_i\alpha_i = \sum_{i=1}^n n_i l_i$, where $\sum_{i=1}^n n_i =$ even.
There are three maximal subgroups of the group $W(d_n)$, two symmetric groups $W(a_{n-1})$ of order $n!$ and one group $W(d_{n-1})$ of order $2^{n-2}(n-1)!$. The facets of the root polytope are of two types of $(n-1)$-ambo simplex and one facet of cross polytope $\beta_{n-1}$. The fundamental weights orthogonal to two ambo simplexes are $\omega_{n-1}$ and $\omega_n$ and the vector orthogonal to the cross polytope $\beta_{n-1}$ is the fundamental weight $\omega_1$. The hyperplane determined by the vectors

$$W(d_n)\{\omega_1, \omega_{n-1}, \omega_n\} = (\omega_1)_{d_n} \cup (\omega_{n-1})_{d_n} \cup (\omega_n)_{d_n}, \quad (26)$$

is orthogonal to the fundamental weight $\omega_2$ and (26) determines the vertices of the dual polytope, the Voronoi cell $V(0)$. Since

$$(\omega_1)_{d_n} = \{\pm l_1, \pm l_2, \ldots, \pm l_n\},\ (\omega_{n-1})_{d_n} \cup (\omega_n)_{d_n} = \tfrac{1}{2}\{\pm l_1 \pm l_2 \pm \cdots \pm l_n\}, \quad (27)$$

the orbit $(\omega_1)_{d_n} = \beta_n$ is a cross polytope and $(\omega_{n-1})_{d_n} \cup (\omega_n)_{d_n} = h\gamma_n \cup h\gamma_n$ representing the union of two hemicubes is actually a cube. Then the volume of the Voronoi cell is

$$Vol\ V(0) = Vol\beta_n + Volh\gamma_n + Volh\gamma_n = 2 \quad (28)$$

as can be checked from Appendix A.

We wonder whether there exits any other fundamental weight in the same hyperplane. When we check $[(\lambda\omega_i - \omega_1), \omega_2] = 0$, $(i = 2, 3, \ldots, n-2)$ we obtain that $\lambda = \frac{1}{2}$. Therefore the vectors $\frac{\omega_2}{2}, \frac{\omega_3}{2}, \ldots, \frac{\omega_{n-2}}{2}$ are in the hyperplane of the Voronoi cell but they do not constitute the vertices of the Voronoi cell. As we will discuss they represent the vertices of the fundamental simplex.

Now, we will discuss the detailed structures of the Voronoi and Delone cells for some simple cases $D_3$ and $D_4$.

i) $D_3 \approx A_3$

The Coxeter-Dynkin diagram is shown in Fig. 5.

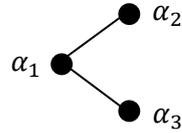

**Figure 5**
Coxeter-Dynkin diagram of $d_3 \approx a_3$.

The root lattice $D_3$ is the face centered cubic (f.c.c.) lattice which has many applications in condensed matter physics and chemistry. As it is well known its Voronoi cell (Wigner-Seitz cell) is the rhombic dodecahedron. The vertices of the Wigner-Seitz cell follows from the union of the orbits

$$V(0) = (\omega_1)_{d_3} \cup (\omega_2)_{d_3} \cup (\omega_3)_{d_3} = \{\pm l_1, \pm l_2, \pm l_3\} \cup \{\tfrac{1}{2}(\pm l_1 \pm l_2 \pm l_3)\}. \tag{29}$$

A plot of the rhombic dodecahedron is shown in Fig. 6.

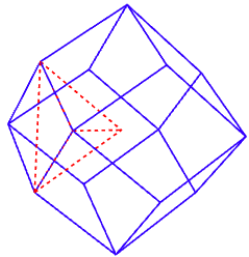

**Figure 6**
Rhombic dodecahedron (the Voronoi cell of f.c.c. lattice) with fundamental simplex identified with the dashed lines.

The Delone cells centered at the origin and constituting the vertices of the Voronoi cell represent an octahedron, a tetrahedron and an inverted tetrahedron respectively. Vertices of the Delone cells centered at the vertices of the Voronoi cell $V(0)$ are obtained by adding the vertices of the orbits $(\omega_1)_{d_3} + (\omega_1)_{d_3}$ to obtain 6 octahedra and $(\omega_2)_{d_3} + (\omega_3)_{d_3}$ to obtain 8 tetrahedra tiling the root lattice closest to the Voronoi cell $V(0)$. Tessellation of the lattice $D_3$ with Delone cells is obtained by the affine group $W(D_3)$. Delone cells centered at the vertices of the rhombic dodecahedron is depicted in Fig.7.

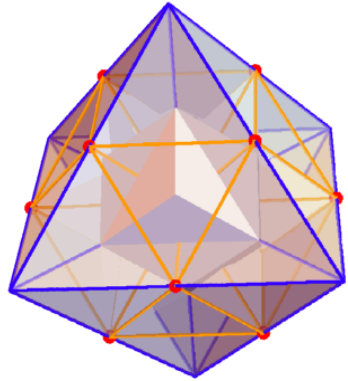

**Figure 7**
Delone cells centered at the Vertices of the rhombic dodecahedron $V(0)$ (the Voronoi cell of f.c.c. lattice centered at the origin), 8 tetrahedra and 6 octahedra surrounds the Wigner-Seitz cell of the f.c.c. lattice.

Since the root polytope (cuboctahedron) is obtained from the highest weight $\omega_2 + \omega_3$ and has 12 vertices the plane orthogonal to this vector is generated by $<r_1>(\omega_1, \omega_2, \omega_3) = \{\omega_1, r_1\omega_1, \omega_2, \omega_3\}$. These 4 vertices determine a rhombus of edge length $\frac{\sqrt{3}}{2}$ with the diagonals of lengths 1 and $\sqrt{2}$ and an area of $\frac{1}{\sqrt{2}}$ as shown in Fig.8. The volume of the pyramid with the rhombic base equals $\frac{1}{3}\frac{\|\omega_1+\omega_2\|}{2}\frac{1}{\sqrt{2}} = \frac{1}{6}$. Since we have 12 pyramid constituting the Voronoi cell the volume equals $Vol\ V(0) = 2$.

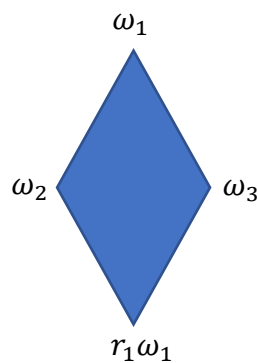

**Figure 8**
One of the facet of the Wigner-Seitz cell with vertices ($\omega_1, \omega_2, \omega_3$ and $r_1\omega_1$).

The fundamental simplex is the semiregular tetrahedron consisting of the vertices $(0, \omega_1, \omega_2, \omega_3)$ with isosceles triangular faces of edge lengths $(\frac{\sqrt{3}}{2}, \frac{\sqrt{3}}{2}, 1)$. Note also that the vector $\frac{1}{2}(l_1 + l_2)$ represents the mid-point of the edge joining $\omega_2$ and $\omega_3$. The roof of the fundamental simplex is the triangle of edge length $(\frac{\sqrt{3}}{2}, \frac{\sqrt{3}}{2}, 1)$, half the rhombus. It is clear that the volume of the fundamental simplex equals $Vol$ (fundamental simplex) $= \frac{1}{12}$ so that $Vol\ V(0) = Vol[W(d_3)(0, \omega_1, \omega_2, \omega_3)] = \frac{24}{12} = 2$.

ii) $D_4 \approx D_4{}^*$

It is a self-dual lattice as will be explained shortly with a point symmetry $W(d_4)$ of order 192. The Coxeter-Dynkin diagram is shown in Fig. 9.

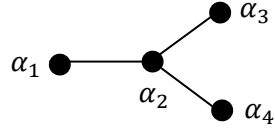

**Figure 9**
Coxeter-Dynkin diagram of $d_4$.

Note that the Dynkin diagram symmetry of $d_4$ leads to the automorphism group $Aut(d_4) \approx W(d_4){:}S_3 \approx W(f_4)$ of order $192 \times 6 = 1152$. The fundamental weights are given by

$$\omega_1 = l_1,\ \omega_2 = l_1 + l_2,\ \omega_3 = \frac{1}{2}(l_1 + l_2 + l_3 - l_4),\ \omega_4 = \frac{1}{2}(l_1 + l_2 + l_3 + l_4). \tag{30}$$

The root polytope is the 24-cell with 24 octahedral facets and 24 vertices represented by the roots $\pm l_i \pm l_j$, ($i \neq j = 1,2,3,4$) of edge length $\sqrt{2}$. Number of its facets equal $N_0 = 24$, $N_1 = 96$, $N_2 = 96$, $N_3 = 24$ satisfying the Euler characteristic equation.

The Voronoi cell is the union of three orbits $V(0) = (\omega_1)_{d_4} \cup (\omega_3)_{d_4} \cup (\omega_4)_{d_4}$ each of which is identical to a cross polytope $\beta_4$. Although the last two is known as hemicubes $h\gamma_4$ but $\beta_4$ and $h\gamma_4$ are congruent polytopes. Volume of the Voronoi cell is three times the volume of the cross polytope $\beta_4$ (see appendix A) that is, $Vol\ V(0) = 3 \times \frac{2}{3} = 2$ which can also be obtained as the determinant of the generator matrix

$$\begin{bmatrix} 1 & -1 & 0 & 0 \\ 0 & 1 & -1 & 0 \\ 0 & 0 & 1 & -1 \\ 0 & 0 & 1 & 1 \end{bmatrix}. \tag{31}$$

The Voronoi cell consists of 24 vertices

$$\{\pm l_1, \pm l_2, \pm l_3, \pm l_4\},\ \frac{1}{2}\{\pm l_1 \pm l_2 \pm l_3 \pm l_4\}, \tag{32}$$

representing another 24-cell of edge length 1, dual to the root polytope. Hence the 24-cell is a self-dual polytope. When the vectors in (32) are represented by quaternions [Koca, Koc & Al-Barwani, 2006] they describe the binary tetrahedral group and together with the normalized set of quaternionic root system of the root polytope they describe the binary octahedral group of quaternions of order 48. Vertices of the facet of the Voronoi cell orthogonal to the highest weight vector $\omega_2$ are obtained by the orbit $< r_1, r_3, r_4 >$ ($\omega_1, \omega_3, \omega_4$) generating 6 vectors

$$l_1, l_2, \frac{1}{2}\{l_1 + l_2 \pm l_3 \pm l_4\} \tag{33}$$

which represents an octahedron centered at $\frac{\omega_2}{2}$ as stated earlier. Its edge length is 1 and its volume is $\frac{\sqrt{2}}{3}$ so that the four-dimensional pyramid based on this octahedron has the volume $\frac{1}{4} \frac{\|\omega_2\|}{2} \frac{\sqrt{2}}{3} = \frac{1}{12}$ and since the Voronoi cell contains 24 such pyramids, the $Vol\, V(0) = 2$. The Delone cells centered at the vertices of the Voronoi polytope $V(0)$ can be determined by the same simple technique, namely, by adding the vertices of $(\omega_1)_{d_4} + (\omega_1)_{d_4}$, $(\omega_3)_{d_4} + (\omega_3)_{d_4}$ and $(\omega_4)_{d_4} + (\omega_4)_{d_4}$. Just to give two examples the vertices $0, 2l_1, l_1 \pm l_2, l_1 \pm l_3, l_1 \pm l_4$ represent the 4-octahedron $\beta_4$ centered at $l_1$ which is obtained by adding $l_1$ on $(\omega_1)_{d_4}$ and the vectors $0, l_1 + l_2, l_1 + l_3, l_1 + l_4, l_2 + l_3, l_2 + l_4, l_3 + l_4, l_1 + l_2 + l_3 + l_4$ are the vertices of the 4-octahedron centered at $\frac{1}{2}(l_1 + l_2 + l_3 + l_4)$ and obtained by adding $\frac{1}{2}(l_1 + l_2 + l_3 + l_4)$ on $(\omega_4)_{d_4}$. This shows how the Voronoi cell $V(0)$ is surrounded by 24 4-octahedra.

The weight lattice is represented by the vector $\sum_{i=1}^{4} c_i \omega_i = \sum_{i=1}^{4} n_i l_i$ with $\sum_{i=1}^{4} n_i =$ (even or odd). Multiplying the vectors by 2 $\sum_{i=1}^{4} n_i =$ (even) we obtain the weight lattice congruent to the root lattice.

The fundamental simplex has the vertices $0, \omega_1 = l_1, \frac{\omega_2}{2} = \frac{l_1 + l_2}{2}, \omega_3 = \frac{1}{2}(l_1 + l_2 + l_3 - l_4)$ and $\omega_4 = \frac{1}{2}(l_1 + l_2 + l_3 + l_4)$ where $\frac{\omega_2}{2}$ represents the center of the octahedron. The roof of the fundamental simplex is a prism with a right triangular base dividing octahedron into 8 congruent pieces. This number actually equals the order of the group $|< r_1, r_3, r_4 >| = 8$. Therefore the volume of the fundamental simplex is $\frac{1}{8} \times \frac{1}{12} = \frac{1}{96}$ which verifies the equation that the facet of the Voronoi cell has the vertices $V(0) = W(d_4)(\omega_1, \omega_3, \omega_4)$ and $Vol\, V(0) = 2$. As we see clearly that the facet of the Voronoi cell of $D_4$ which is an octahedron changed radically as confronted to the rhombic face of the Voronoi cell of $A_4$. However both facets have similarities with $D_3 \approx A_3$. In the case of $A_4$ the rhombus of $A_3$ generalizes to rhombohedron. We can regard the rhombus of Fig.8 as dipyramid with a basis of line segment and it generalizes to a dipyramid (octahedron) with a square base as one moves from the line segment to a square by increasing the dimension one more. As we will see in the most general case of the root lattice $D_n$, the facet of the Voronoi cell is a dipyramid with a base of $(n-2)$-cube.

Let us have a look at another simpler case $D_5$ before we discuss the general case. The vertices of the facet of the lattice $D_5$ orthogonal to the weight vector $\omega_2$ are generated by $< r_1, r_3, r_4, r_5 >$ $(\omega_1, \omega_4, \omega_5)$ leading to the vertices

$$l_1, l_2, \frac{1}{2}\{l_1 + l_2 \pm l_3 \pm l_4 \pm l_5\}. \tag{34}$$

This is a dipyramid with a base of cube with isosceles triangular 2d-faces of edge lengths $(\frac{\sqrt{5}}{2}, \frac{\sqrt{5}}{2}, 1)$ and its volume is $\frac{\sqrt{2}}{4}$. The fundamental simplex has the vertices $0, l_1, \frac{1}{2}(l_1 + l_2), \frac{1}{2}(l_1 + l_2 + l_3), \frac{1}{2}(l_1 + l_2 + l_3 + l_4 \pm l_5)$ where two vertices $\frac{\omega_2}{2} = \frac{l_1 + l_2}{2}$ and $\frac{\omega_3}{2} = \frac{1}{2}(l_1 + l_2 + l_3)$ correspond to the center and one of the face of the cube respectively and they are not the vertices of the facet. The isosceles triangles of vertices (34) project onto the Coxeter plane as two triangles one with isosceles triangle with interior angles $\frac{\pi}{8}, \frac{6\pi}{8}$ and edge lengths $\sin\left(\frac{\pi}{8}\right)$ and $\sin\left(\frac{2\pi}{8}\right)$. The second triangle has interior angles $\frac{\pi}{8}, \frac{2\pi}{8}, \frac{5\pi}{8}$ and the corresponding edge lengths are $\sin\left(\frac{\pi}{8}\right)$, $\sin\left(\frac{2\pi}{8}\right)$, $\sin\left(\frac{3\pi}{8}\right)$. Some examples of the aperiodic tilings by these two tiles are depicted in Fig.10.

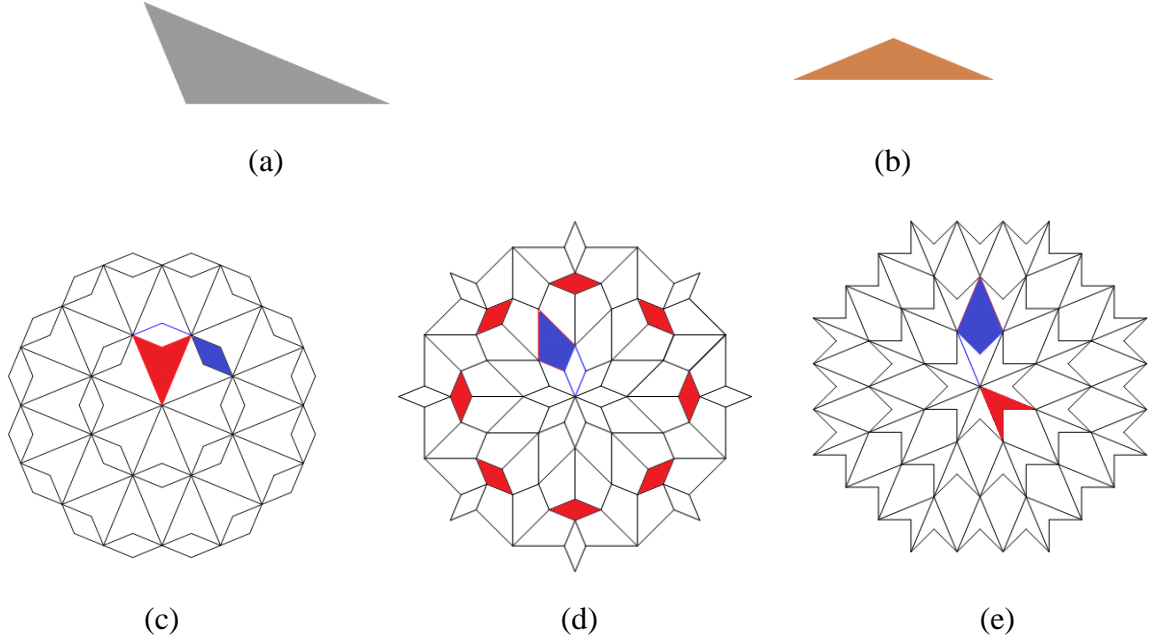

(a)            (b)

(c)            (d)            (e)

**Figure 10**
Tiles (a, b) and examples of 8-fold symmetric aperiodic tilings by projection of $D_5$ root lattice (c, d, e). (All tiles in c, d, e are made of the basic tiles a and b).

Now we can briefly discuss the Voronoi and Delone cells of the lattice $D_n$ in a general context. The number of facets of the root polytope can be determined as

$$N_0^r = \frac{2^{n-1} n!}{2^{n-2}(n-2)!} = 2n(n-1),$$

$$N_1^r = \frac{2^{n-1} n!}{2^{n-3}(n-3)!} = 2^2 n(n-1)(n-2),$$

$$N_d^r = 2^{d+1} \binom{n}{n-d-1}(2(n-d-1)+1), \ 2 \leq d \leq n-3, \tag{35}$$

$$N_{n-2}^r = 3 \times 2^{n-1} n,$$

$$N_{n-1}^r = 2^n + 2n.$$

Since the number of vertices of the facets of the Voronoi polytope is given by $N_d^v = N_{(n-d-1)}^r$ we obtain the vertices of the Voronoi cell as

$$N_0^v = 2^n + 2n,$$

$$N_1^v = 3 \times 2^{n-1} n,$$

$$N_d^v = 2^{n-d} \binom{n}{d}(2d+1), \ 2 \leq d \leq n-3, \tag{36}$$

$$N_{n-2}^v = 2^2 n(n-1)(n-2),$$

$$N_{n-1}^v = 2n(n-1).$$

Structure and vertices of the Voronoi cell were already given in (26-27). To determine the Delone cells centered at the vertices of the Voronoi cell $V(0)$ we repeat what we did for $D_3$ and $D_4$. So there are two different cases:

*i)* even $n$

Vertices of the Delone cells centered at the origin are added to themselves as

$$(\omega_1)_{d_n} + (\omega_1)_{d_n}, (\omega_{n-1})_{d_n} + (\omega_{n-1})_{d_n}, (\omega_n)_{d_n} + (\omega_n)_{d_n}. \tag{37}$$

As an example, the vertices of the Delone cell centered at the vertex $\omega_n = \frac{1}{2}(l_1 + l_2 + \cdots + l_{n-1} + l_n)$ are given by

$$0, l_{i_1} + l_{i_2}, l_{i_1} + l_{i_2} + l_{i_3} + l_{i_4}, \ldots, l_{i_1} + l_{i_2} + \cdots + l_{i_{n-2}},$$
$$l_1 + l_2 + \cdots + l_n, \; i_1 \neq i_2 \neq \cdots \neq i_{n-2} \tag{38}$$

where the number of vertices in (38) are

$$\binom{n}{0} + \binom{n}{2} + \binom{n}{4} + \cdots + \binom{n}{n-2} + \binom{n}{n} = 2^{n-1}. \tag{39}$$

*ii)* odd $n$

Delone cells centered at the origin are added as $(\omega_1)_{d_n} + (\omega_1)_{d_n}, (\omega_{n-1})_{d_n} + (\omega_n)_{d_n}$.
A similar set of vertices for the Delone cell centered at $\omega_n = \frac{1}{2}(l_1 + l_2 + \cdots + l_{n-1} + l_n)$ for odd $n$ is obtained where the number of vertices are calculated as

$$\binom{n}{0} + \binom{n}{2} + \binom{n}{4} + \cdots + \binom{n}{n-3} + \binom{n}{n-1} = 2^{n-1}. \tag{40}$$

We now return back to the structure of the $(n-1)$-facet of the Voronoi cell $V(0)$. Vertices of the facet orthogonal to the vector $\omega_2$ is obtained from

$$< r_1, r_3, r_4, \ldots, r_{n-1}, r_n > (\omega_1, \omega_{n-1}, \omega_n) = \{\omega_1, r_1 \omega_1, < r_3, r_4, \ldots, r_{n-1}, r_n > (\omega_{n-1}, \omega_n)\}. \tag{41}$$

The vertices, $< r_3, r_4, \ldots, r_{n-1}, r_n > (\omega_{n-1}, \omega_n)$ constitute the union of two hemicubes equivalent to a cube in $(n-2)$ dimensions with vertices $\frac{1}{2}\{l_1 + l_2 \pm l_3 \pm \cdots \pm l_{n-1} \pm l_n\}$. The vectors $\omega_1 = l_1, r_1 \omega_1 = l_2$ represent two opposite vertices of dipyramid where the center of the cube is given by $\frac{\omega_2}{2} = \frac{l_1+l_2}{2}$. The volume of this dipyramid equals $2 \times \frac{1}{n-1} \times \frac{1}{\sqrt{2}}(1)^{n-2} = \frac{\sqrt{2}}{n-1}$. Since the Voronoi cell consists of $2n(n-1)$ pyramids with the $(n-1)$-dimensional dipyramid as a base and height $\frac{\|\omega_2\|}{2} = \frac{1}{\sqrt{2}}$ then the volume of the Voronoi cell is given as $\frac{2n(n-1)}{n} \times \frac{1}{\sqrt{2}} \frac{\sqrt{2}}{n-1} = 2$ as expected. The vertices $\frac{\omega_2}{2}, \frac{\omega_3}{2}, \ldots, \frac{\omega_{n-2}}{2}$ of the fundamental simplex are not the vertices of the Voronoi cell rather they correspond to the center, and centers of faces of the $(n-2)$-dimensional cube. Volume of the roof of the fundamental simplex is given by $\frac{\sqrt{2}}{2^{n-2}(n-1)!}$. Multiplying this by the height $\frac{1}{\sqrt{2}}$ and dividing by $n$ gives the volume of the fundamental simplex, $\frac{1}{2^{n-2}n!}$.

## 5. The Weight Lattice $A_n^*$ and related polytopes

The numbers of facets of the weight lattice $A_n^*$ have been studied by Louis Michel in details in [Michel, 1997]. The weight lattice is represented by the vectors as linear combinations of the fundamental weights $q \equiv (q_1, q_2, \ldots, q_{n+1}) = \sum_{j=1}^{n} b_j \omega_j$, $b_j \in \mathbb{Z}$ where $q_i$ are the components in the $l_i$-bases. As can be seen from (8), each component is given by $q_i = \frac{m_i}{n+1}$, $m_i \in \mathbb{Z}$ and satisfies $(q_i - q_j) \in \mathbb{Z}$. The center of the fundamental simplex defined by $Q = \frac{(0+\omega_1+\omega_2+\cdots+\omega_n)}{n+1}$ is at its maximal distance from the vertices of the fundamental simplex, that is, it is one of the holes of the weight lattice. Therefore the center of the fundamental simplex represents one of the vertex of the Voronoi polytope of the weight lattice. The components of $Q$ in the $l_i$-basis reads

$$Q = \frac{1}{n+1}\left[\frac{n}{2}, \left(\frac{n}{2}-1\right), \left(\frac{n}{2}-2\right), \ldots, \left(1-\frac{n}{2}\right), -\frac{n}{2}\right]. \tag{42}$$

The orbit $W(a_n)Q$ represents the $(n+1)!$ vertices of the Voronoi cell centered at the origin $V(0)^* = W(a_n)Q$. The Coxeter-Weyl group $W(a_n)$ permutes the $n+1$ components of $Q$ and for this reason the polytope is also called the permutohedron. For further discussions of the Voronoi cell of $A_n^*$ we refer the reader to the references [Vallentin, 2003] and [Garber, 2012]. We recall that because of the Dynkin diagram symmetry the lengths of the fundamental weights pairwise equal each other

$$\|\omega_1\| = \|\omega_n\|, \; \|\omega_2\| = \|\omega_{n-1}\|, \ldots \tag{43}$$

and satisfy the inequality

$$\|\omega_1\| < \|\omega_2\| < \|\omega_3\| < \cdots. \tag{44}$$

To elaborate the topic some examples are in order.

i) $A_3^*$ lattice

Since $A_3 \approx D_3$ we will use three-dimensional representation of the roots and fundamental weights defined by

$$\alpha_1 = l_2 - l_3, \alpha_2 = l_1 - l_2, \alpha_3 = l_2 + l_3,$$
$$\omega_1 = \frac{1}{2}(l_1 + l_2 - l_3), \omega_2 = l_1, \omega_3 = \frac{1}{2}(l_1 + l_2 + l_3). \tag{45}$$

The center of the fundamental simplex now reads $Q = \frac{1}{4}(\omega_1 + \omega_2 + \omega_3) = \frac{1}{4}(2l_1 + l_2)$.

The generators of the Coxeter-Weyl group operates like $r_1: l_2 \leftrightarrow l_3$, $r_2: l_1 \leftrightarrow l_2$, $r_3: l_2 \leftrightarrow -l_3$.

Therefore, the orbit is a polytope with 24 vertices given in terms of components

$$V(0)^* = W(a_3)Q$$
$$= \left\{\frac{1}{4}(\pm 2, \pm 1, 0), \frac{1}{4}(\pm 1, \pm 2, 0), \frac{1}{4}(\pm 2, 0, \pm 1), \frac{1}{4}(\pm 1, 0, \pm 2), \frac{1}{4}(0, \pm 2, \pm 1), \frac{1}{4}(0, \pm 1, \pm 2)\right\}. \tag{46}$$

This is the truncated octahedron representing the Wigner-Seitz cell of the b.c.c. lattice. It is vertex transitive but not face transitive as it consists of 8 hexagonal faces in two orbits and 6 square faces as another orbit under the group $W(a_3)$ as shown in Fig.11.

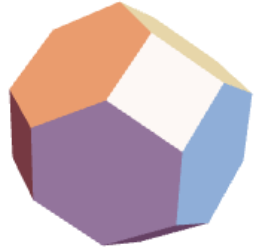

**Figure 11**
Truncated octahedron as the Voronoi cell of the b.c.c. lattice.

The center of the hexagon whose vertices are generated by the subgroup $<r_1, r_2> \frac{1}{4}(\omega_1 + \omega_2 + \omega_3)$ is the vector $\frac{\omega_3}{2}$ of length $\frac{\sqrt{3}}{4} \approx 0.43$ and the center of the square generated by the group $<r_1, r_3> \frac{1}{4}(\omega_1 + \omega_2 + \omega_3)$ is the vector $\frac{\omega_2}{2}$ of length $\frac{1}{2} = 0.5$. Had we operated the subgroup $<r_2, r_3> \frac{1}{4}(\omega_1 + \omega_2 + \omega_3)$ we would have obtained another hexagon with the vector $\frac{\omega_1}{2}$ denoting its center. The 8 hexagons correspond to the orbits of these fundamental weights so that the in-sphere of the Voronoi cell will touch the 8 hexagonal surfaces. The spheres touching the in-sphere has centers represented by vertices $W(a_3)(\omega_1)$ and the $W(a_3)(\omega_3) = -W(a_3)(\omega_1)$. They are representing two tetrahedra forming a cube with the vertices $\frac{1}{2}(\pm 1, \pm 1, \pm 1)$. Such a double simplex is named as diplo-simplex [Conway & Sloane, 1991] as we will see in a more general case of $A_n^*$. The diplo-simplex $(\omega_1)_{a_3} \cup (\omega_3)_{a_3}$ is the contact polytope of the $A_3^*$ lattice. By applying the group $W(a_3)$ on the fundamental simplex one obtains 24 copies of it, each of which is centered at one of the vertex of the Voronoi cell implying that the Delone cell of the weight lattice is the fundamental simplex.

*ii*) $A_4^*$ lattice

For $n \geq 4$ we use the usual representations with $n + 1$ orthonormal vectors for convenience although it is always possible to find a basis with $n$ orthonormal vectors. As we have seen in Sec.3 they are not as practical as $n + 1$ coordinate system except only for the case of $A_4^*$ one can use quaternions to describe the lattice as well as the point group [Koca, Ozdes Koca & Al-Ajmi, 2012]. In this reference all polytopes symmetric under the Coxeter-Weyl group $W(a_4)$ have been studied in details. The permutohedron $\frac{(1111)_{a_4}}{5}$ has facets of the form of truncated octahedra and hexagonal prisms. The centers of the truncated octahedra are represented by the union of the orbits $\frac{1}{2}[(\omega_1)_{a_4} \cup (\omega_4)_{a_4}]$ and the centers of the hexagonal prisms by the union of the orbits $\frac{1}{2}[(\omega_2)_{a_4} \cup (\omega_3)_{a_4}]$ and because of (42-43) the 10 spheres touching the in-sphere has the orbit of the diplo simplex $(\omega_1)_{a_4} \cup (-\omega_1)_{a_4}$. A few words are in order why $\frac{\omega_4}{2}$ represents the center of the truncated octahedron $<r_1, r_2, r_3> \frac{1}{5}(\omega_1 + \omega_2 + \omega_3 + \omega_4)$. Average of these 24 vectors is proportional to $\omega_4$ since it is invariant under the group $<r_1, r_2, r_3>$. Then the center of the truncated octahedron can be written as

$$c = \frac{1}{5}[(\omega_1 + \omega_2 + \omega_3 + \omega_4) + m_1\alpha_1 + m_2\alpha_2 + m_3\alpha_3)] \tag{47}$$

Using (2) to eliminate the vectors $\omega_1, \omega_2,$ and $\omega_3$, one determines the coefficients as $m_1 = m_3 = -\frac{3}{2}$, $m_2 = -2$ and finds $c = \frac{\omega_4}{2}$ after substitution. This makes it clear that the contact polytope is the diplo-simplex

$$(\omega_1)_{a_4} \cup (-\omega_1)_{a_4} = \{\pm\omega_1, \pm(\omega_1 - \omega_2), \pm(\omega_2 - \omega_3), \pm(\omega_3 - \omega_4)\}. \tag{48}$$

The Delone cell of $A_4^*$ is the fundamental simplex with vertices $(0, \omega_1, \omega_2, \omega_3, \omega_4)$ which is in turn a 4-simplex of identical irregular tetrahedral facets where identical tetrahedra consist of isosceles triangles of two types. Volume of the permutohedron can also be calculated as $\frac{1}{\sqrt{5}}$ by using the volumes of its cells.

It is clear that the polytopes of the $A_n^*, (n > 4)$ follows the same pattern of logic. The list of numbers of the facets of the Voronoi cell (permutohedron) of $A_n^*$ can be found in [Michel, 1995]. For example, the $(n-1)$ facets generated by

$$< r_1, r_2, \ldots, r_{n-1} > \frac{1}{n+1}(\omega_1 + \omega_2 + \cdots + \omega_{n-1} + \omega_n) \tag{49}$$

or by its conjugate groups are the permutohedra and their centers can be determined as the orbits of the vectors $\frac{\omega_n}{2}$ or $\frac{\omega_1}{2}$. Centers of the other facets of the Voronoi cell equals the orbits $(\frac{\omega_i}{2})_{a_n}, i \neq 1, n$ and the their lengths are greater than the lengths of $\frac{\omega_n}{2}$ and $\frac{\omega_1}{2}$ as shown in (43-44). Therefore the contact polytope is the diplo simplex $(\omega_1)_{a_n} \cup (-\omega_1)_{a_n}$. The proof goes as claimed in (47) where the center of one of the $(n-1)$-permutohedra is given by

$$c = \frac{1}{n+1}[(\omega_1 + \omega_2 + \cdots + \omega_{n-1} + \omega_n) + m_1\alpha_1 + \cdots + m_{n-2}\alpha_{n-2} + m_{n-1}\alpha_{n-1}]$$
$$= \frac{1}{n+1}(1 - m_{n-1})\omega_n. \tag{50}$$

Since (50) is invariant under the group $< r_1, r_2, \ldots, r_{n-1} >$ it is independent of all fundamental weights except $\omega_n$ which leads to $(n-1)$ linear equations in $m_1, \ldots, m_{n-1}$ and a solution is obtained as $m_{n-1} = \frac{1-n}{2}$ yielding $c = \frac{\omega_n}{2}$.

The summary of the polytopes of $A_n^*$ lattice can be stated as follows.

The Delone cell is the fundamental simplex : $(0, \omega_1, \omega_2, \ldots, \omega_{n-1}, \omega_n)$;
The Voronoi cell is the permutohedron (orbit of the scaled Weyl vector) with vertices: $\frac{(11\ldots11)_{a_n}}{n+1}$;
The Contact polytope is the diplo simplex: $\{\pm\omega_1, \pm(\omega_1 - \omega_2), \ldots, \pm(\omega_{n-1} - \omega_n), \pm\omega_n\}$.

## 6. The Weight Lattice $D_n^*$ and related polytopes

Since $n \leq 4$ is already studied in Sec.4 our discussion in this section is valid for $n \geq 5$. This section directly follows the paper [Conway & Soane, 1991] with some additional remarks and examples such as $D_5^*$ and $D_6^*$. To determine the vertices of the Voronoi cell $V(0)^*$ we first determine one of its vertex as the vector $P$ equidistant from the four vertices $(0, \omega_1, \omega_{n-1}, \omega_n)$. If $n = 2t$ even, the vector equidistant from four vertices is $P = \frac{\omega_t}{2}$, which is one of the vertices of the fundamental simplex and the distance to each vector equals $\frac{\sqrt{t}}{2}$. For odd $n = 2t + 1$, $P = \frac{\omega_t + \omega_{t+1}}{4}$ is the point midway between two vertices of the fundamental simplex and its distance to four points is $\frac{\sqrt{4t+1}}{4}$. In either case the vertices of the Voronoi cell is obtained as the orbit $V(0)^* = W(d_n)P$. The facets of the Voronoi cell are of two types to which the fundamental weights $\omega_1, \omega_{n-1}, \omega_n$ are orthogonal. The facets orthogonal to $\omega_{n-1}, \omega_n$ are of the same type. The centers of the facets are at distances $\frac{1}{2} \| \omega_1 \|, \frac{1}{2} \| \omega_{n-1} \| = \frac{1}{2} \| \omega_n \|$ and since
$\| \omega_1 \| < \| \omega_{n-1} \| = \| \omega_n \|$, the contact polytope is the $n$-octahedron with $2n$ vertices

$$(\omega_1)_{d_n} = \{\pm l_1, \pm l_2, \dots, \pm l_n\}. \tag{51}$$

The Delone cells are the lattice points nearest $P$. If $n = 2t$, the vertices of the Delone cell centered at $P$ are the vectors given by

$$P + \frac{1}{2}(\pm l_1 \pm l_2 \pm \cdots \pm l_t), \quad P + \frac{1}{2}(\pm l_{t+1} \pm l_{t+2} \pm \cdots \pm l_{2t}). \tag{52}$$

They represent two hypercubes in complementary $t$-spaces sharing the same center $P$ and are called *join* of two hypercubes. If $n = 2t + 1$, the vertices of the Delone cell are given by

$$P + \frac{1}{2}(\pm l_1 \pm l_2 \pm \cdots \pm l_t), \quad P + \frac{1}{4}l_{t+1} + \frac{1}{2}(\pm l_{t+2} \pm l_{t+3} \pm \cdots \pm l_{2t+1}). \tag{53}$$

They represent two hypercubes in orthogonal $t$-spaces whose centers are separated by the vector $\frac{1}{4}l_{t+1}$ and they are called *separated join* of two hypercubes. Some examples are given as follows.

i) $D_5^*$ lattice
The Voronoi cell is generated by the vector $P = \frac{\omega_2 + \omega_3}{4} = \frac{1}{2}(l_1 + l_2 + \frac{1}{2}l_3)$ which is left invariant by the subgroup $<r_1, r_4, r_5>$ of order 8. Therefore the number of vertices are given by $\frac{2^4 5!}{2^3} = 240$. It consists of two types of $a_4$ polytopes and one type of $d_4$ polytope as 4-facets. Its contact polytope is the 5-octahedron with 10 vertices. The Delone cell centered at $P$ is a *separated join* of two squares in two complementary spaces whose vertices are given by

$$P + \frac{1}{2}(\pm l_1 \pm l_2), P + \frac{1}{4}l_3 + \frac{1}{2}(\pm l_4 \pm l_5). \tag{54}$$

ii) The lattice $D_6^*$

The Coxeter-Weyl group $W(d_6)$ admits the icosahedral group as a maximal subgroup and its lattices are important from the point of view of icosahedral quasicrystallography [Koca, N. O., Koca, M.& Koc, R., 2015]. Its Voronoi cell is a polytope whose vertices are generated by the group $W(d_6)\frac{\omega_3}{2}$ where $P = \frac{1}{2}(l_1 + l_2 + l_3)$. It has 160 vertices and consists of second order diplo simplexes and the Voronoi polytope of $D_5$ as 5-facets. The contact polytope is a 6-octahedron with 12 vertices and it represents an icosahedron when it is projected into 3-space. The Delone cell is the join of two cubes in two complementary spaces whose vertices are given by

$$P + \frac{1}{2}(\pm l_1 \pm l_2 \pm l_3), P + \frac{1}{2}(\pm l_4 \pm l_5 \pm l_6). \tag{55}$$

This is a cube in 6-dimensions whose center is shifted to the point $P$. Further studies of the projection of the lattices $D_6$ and $D_6^*$ and their polytopes under the icosahedral symmetry could be interesting for the quasicrystallographic icosahedral crystals.

## 7. Concluding Remarks

We have presented a detailed study of the root and weight lattices of the $A$-$D$ series by highlighting some examples which could be useful in quasicrystallography. The special polytopes such as Voronoi cells, contact polytopes (root polytope in the case of root lattice), Delone cells have been identified. The explicit structures of the facets of the Voronoi and Delone cells have been worked out. For the first time it was noted that the facets of the Voronoi cell of the lattice $A_n$ are the generalized rhombohedra and those of lattice $D_n$ are the dipyramids based on the hypercubes. Tessellation by Delone cells have been exemplified in many cases. Volumes of the Voronoi cells are calculated via their facets and volumes of certain regular polytopes such as $\alpha_n, \beta_n$ and $h\gamma_n$ have been calculated, some with via recurrence relations.

## Appendix A: Volumes of certain polytopes

1) $n$-simplex $((n+1)$-cell$)$ $\alpha_n$

It is the orbit $(\omega_1)_{a_n}$ or $(\omega_n)_{a_n}$. If we take $(\omega_1)_{a_n}$ as the $n$-simplex its vertices are given by

$$(\omega_1)_{a_n} = \{\omega_1, \omega_2 - \omega_1, \omega_3 - \omega_2, \ldots, \omega_n - \omega_{n-1}, -\omega_n\}. \tag{A1}$$

The simplex consists of $(n+1)$ $\alpha_{n-1}$ facets. Without any loss of generality the first $n$-vertices can be taken as the vertices of the simplex $\alpha_{n-1}$. The volume of $\alpha_n$ is the sum of the volumes of the $(n+1)$ pyramids based on the $\alpha_{n-1}$ facets. The height of the $\alpha_{n-1}$ facet is the average of the first $n$ vertices $h = \frac{\|\omega_n\|}{n} = \frac{1}{n}\sqrt{\frac{n}{n+1}}$. Then the volume of $\alpha_n$ can be written as a recurrence relation,

$$Vol(\alpha_n) = \frac{n+1}{n^2}\sqrt{\frac{n}{n+1}} Vol(\alpha_{n-1}) = \frac{1}{n}\sqrt{\frac{n+1}{n}} Vol(\alpha_{n-1}),$$

$$Vol(\alpha_{n-1}) = \frac{1}{n-1}\sqrt{\frac{n}{n-1}} Vol(\alpha_{n-2}),$$

.
.
. \hfill (A2)
.

$$Vol(\alpha_2) = \frac{1}{2}\sqrt{\frac{3}{2}} Vol(\alpha_1) = \frac{\sqrt{3}}{2},$$

that leads to the result $Vol(\alpha_n) = \frac{\sqrt{n+1}}{n!}$.

2) $n$-cross polytope ($n$-octahedron) $\beta_n$

The $n$-octahedron $\beta_n = (\omega_1)_{d_n}$ has two $\alpha_{n-1}$ facets each of which occurs as many as $\frac{2^{n-1}n!}{n!} = 2^{n-1}$. Therefore $\beta_n$ consists of $2^n \alpha_{n-1}$ cells. The vertices of $\alpha_{n-1}$ can be obtained applying the group element $a = r_1 r_2 \ldots r_{n-1}, (a^n = 1)$ on $\omega_1 = l_1$ $n$ times which will lead to the vertices $l_1, l_2, ..., l_n$. The average of these vectors gives the height $h = \frac{1}{\sqrt{n}}$ of the cell $\alpha_{n-1}$. Therefore the volume equals

$$Vol(\beta_n) = \frac{2^n}{n!}. \tag{A3}$$

3) $n$-hemicube $h\gamma_n$

It is represented by the either orbit $(\omega_{n-1})_{d_n}$ or $(\omega_n)_{d_n}$ and has two types of $(n-1)$- facets; one facet is a simplex $\alpha_{n-1}$ and the other one is a hemicube $h\gamma_{n-1}$. The volume of hemicube then can be written as

$$Vol(h\gamma_n) = \frac{1}{n}\left[N_{\alpha_{n-1}} h_{\alpha_{n-1}} Vol(\alpha_{n-1}) + N_{h\gamma_{n-1}} h_{h\gamma_{n-1}} Vol(h\gamma_{n-1})\right]. \tag{A4}$$

Vertices of $\alpha_{n-1}$ is obtained by applying the group element $a = r_1 r_2 ... r_{n-1}$, on $\omega_{n-1}$, $n$ times and one averages over the vertices to find the center as $c = \frac{1}{n}(n-2)\omega_n$ which leads to the height $h_{\alpha_{n-1}} = \frac{n-2}{2\sqrt{n}}$.

Vertices of $h\gamma_{n-1}$ can be obtained from $< r_2, r_3, ..., r_{n-1}, r_n > \omega_{n-1} = \frac{1}{2}(l_1 \pm l_2 \pm, ..., \pm l_n)$ with odd number of $(-)$ signs. The center of $h\gamma_{n-1}$ is then $c = \frac{1}{2}$. Using (A2) we can write the volume of $h\gamma_n$ as

$$Vol(h\gamma_n) = \frac{1}{n}\left[\frac{2^{n-2}(n-2)}{(n-1)!} + 2n \times \frac{1}{2} Vol(h\gamma_{n-1})\right] = Vol(h\gamma_{n-1}) + \frac{2^{n-2}(n-2)}{n!},$$

$$Vol(h\gamma_{n-1}) = Vol(h\gamma_{n-2}) + \frac{2^{n-3}(n-3)}{(n-1)!}, \tag{A5}$$

$$\vdots$$

$$Vol(h\gamma_4) = Vol(h\gamma_3) + \frac{2^2(2)}{4!}.$$

By adding all terms and using $Vol(h\gamma_3) = \frac{2}{3!}$ we obtain

$$Vol(h\gamma_n) = 1 - \frac{2^{n-1}}{n!}, n \geq 3. \tag{A6}$$

**References**

Baake, M., Joseph, D., Kramer, P. & Schlottmann, M. (1990). *J. Phys. A: Math. & Gen.* **23**, L1037-L1041.

Bourbaki, N. (1968). *Groupes et Algèbres de Lie. Chap. IV-VI, Actualités Scientifiques et Industrielles* (Paris: Hermann) **288**; English translation: *Lie Groups and Lie algebras* (Springer 2nd printing, 2008).

Boyle, L.& Steinhardt, P.J. ( 2016 ). Coxeter pairs, Amman patterns and Penrose-like Tilings, arXiv: 1608.08215.

Chen, L., Moody, R. V. & Patera, J. (1998). *Quasicrystals and Discrete Geometry*, edited by J. Patera, pp. 135–178. Fields Institute Monographs, Vol. 10. Providence, RI: American Mathematical Society.

Conway, J. H. & Sloane, N.J.A. (1988). *Sphere Packings, Lattices and Groups*. Springer-Verlag New York Inc.

Conway, J. H. & Sloane, N.J.A. (1991). The cell structures of certain lattices. Miscellanea Mathematica, ed. Hilton, P., Hirzebruch, F. & Remmert, R.( New York: Springer).

Coxeter, H.S.M. (1946). *Integral Cayley Numbers*, DMJ **13**(4), 561-578.

Coxeter, H. S. M. (1973). *Regular Polytopes*, Third Edition, New York: Dover Publications.

de Brujin, N. G. (1981). *Algebraic theory of Penrose's non-periodic tilings of the plane*. Nederl. Akad. Wetensch. Proceedings Ser. A **84** (=Indagationes Math. **43**), 38.

Delaunay, N. B. (1929). *Sur la partition reguilere de l'espace a 4-dimensions, Izv. Akad. Nauk SSSR Otdel. Fiz.-Mat. Nauk*. pp. 79-110 and 145-164.

Delaunay, N. B. (1938). *Geometry of positive quadratic forms , Uspehi Mat. Nauk* **3,** 16-62; **4**, 102-164.

Deza, M. & Laurent, M. (1997). *Geometry of Cuts and Metrics, Algorithms and Combinatorics* **15**, Springer-Verlag, Berlin.

Deza, M. & Grishukhin,V. (2004). Non-rigidity degrees of root lattices and their duals, Geometriae Dedicata **104**, 15-24, Kluwer.

Duneau, M. & Katz, A. (1985). *Phys. Rev. Lett.* **54**, 2688–2691.

Engel, P. (1986). *Geometric Crystallography: an Axiomatic Introduction to Crystallography*. Dordrecht: Springer.

Engel, P., Michel, L. & Senechal, M. (1994). *Lattice Geometry, preprint* IHES /P/04/45.

Fritzsch, H., Gell-Mann, M. & Leutwyler, H. (1973). *Phys. Lett.* B**47**, 365-368.

Fritzsch, H. & Minkowski, P. (1975). *Ann. Phys.* **93**, 193-266.

Garber, A. (2012), *Belt distance between the facets of space filling zonotopes*, *Mathematical Notes,* **92** (3-4), 345-355.

Georgi, H. & Glashow, S. (1974). *Phys. Rev. Lett*. **32**, 438-441.

Glashow, S. (1961). *Nucl. Phys.* **22** (4), 579-588.

Grunbaum, B. (1967). *Convex Polytopes*, Wiley, N.Y.

Gürsey, F., Ramond, P., Sikivie, P. (1976). *Phys. Lett*. B**60**, 177-180.

Humphreys, J. E. (1992). *Reflection Groups and Coxeter Groups*. Cambridge University Press.

Koca, M. & Ozdes, N. (1989). *J. Phys. A Math. Gen*. A**22**, 1469-1493.

Koca, M., Koc, R. & Al-Barwani, M. (2006). *Quaternionic root systems and subgroups of the Aut* $(F_4)$, J. M. Phys. **47**, 043507.

Koca, M. (2007). *Quaternionic and Octonionic Structures of the Exceptional Lie Algebras* in Mathematical Physics, Proceedings of the 12th Regional Conference (edit. Aslam, M. J., Hussain, F., Qadir, A., Riazuddin & Saleem, H.), World Scientific p.25.

Koca, N. O., Koca, M. & Al-Ajmi, M. (2012). *Int. J. Geom. Methods Mod. Phys.* **9**(4), 1250035.

Koca, M., Koca, N. O. & Koc, R. (2014). *Int. J. Geom. Meth. Mod. Phys*. **11**, 1450031.

Koca, N. O., Koca, M. & Koc, R. (2015), *Acta Cryst* A**71**, 175-185.

Koca, N. O., Koca, M. & Al-Siyabi, A. (2018). *Int. J. Geom. Methods Mod. Phys.***15** (4), 1850058.

Michel, L. (1995). *Bravais classes, Voronoi cells, Delone cells*, Symmetry and Structural Properties of Condensed Matter, ed. Lulek, T., Florek, W. & Walcerz, S., Zajaczkowo, p. 279, Singapore: World Scientific.

Michel, L. (1997). *Complete description of the Voronoi cell of the Lie algebra $A_n$ weight lattice. On the bounds for the number of d-faces of the n-dimensional Voronoi cells,* IHES/P/97/53.

Moody, R.V. & Patera, J. (1992). *J. Phys. A: Math. Gen.* **25**, 5089-5134.

Novoselov, K. S., Geim, A. K., Morozov, S. V., Jiang, D., Zhang, Y., Dubonos, S. V., Grigorieva, I. V. & Firsov, A. A. (2004). *Science* **306**, 666.

Salam, A. (1968*). Elementary Particle Physics: Relativistic Groups and Analyticity*, ed. N. Eighth Nobel Symposium, Stockholm, edited by N. Svartholm, p. 367.

Senechal, M. (1995). *Quasicrystals and Geometry* (Cambridge University Press).

Slansky, R. (1981). *Phys. Rep.* **79**, 1-128.

Vallentin, F. (2003). *Sphere coverings and tilings (PhD thesis).* Fakultät für Mathematik, Technischen Universität München, Berlin

Voronoi, G. (1908). *Recherches sur les paralleloedres primitifs I. Proprietes generals des paralleloedres. J. reine angew. Math.***133**, 198-287.

Voronoi, G. (1909). *Recherches sur les paralleloedres primitifs II. Domaines de formes quadratiques corresondant aux different types de Nouvelles applications des paramétres continus `a l` a théorie des formes quadratiques primitifs. J. reine angew. Math.***136**, 67-181.

Weinberg, S. (1967). *Phys. Rev. Lett.* **19**, 1264-1266.